\documentclass{article}

\usepackage{arxiv}
\usepackage{subfigure}
\usepackage{amsmath,amssymb,euscript,yfonts,psfrag,latexsym,dsfont,graphicx}
\usepackage{bbm,color,amstext,wasysym,parskip,balance,cite}
\usepackage{comment}

\graphicspath{{./},{./pictures/}}

\newcommand{\PFO}{Perron-Frobenius operator}
\newcommand{\X}{\mathbb{X}} 
\newcommand{\Rd}{\mathbb{R}^d}
\newcommand{\R}{\mathbb{R}}
\newcommand\PD{\mathcal{P}_2(\mathbb{X})}

\newtheorem{theorem}{Theorem}
\newtheorem{lemma}[theorem]{Lemma}
\newtheorem{proposition}[theorem]{Proposition}
\newtheorem{problem}{Problem}
\newtheorem{remark}{Remark}

\title{Data-Driven Approximation of the Perron-Frobenius Operator Using the Wasserstein Metric}

\author{
  Amirhossein Karimi \\
  Department of Mechanical\\
  and Aerospace Engineering\\
  University of California, Irvine\\
  California, USA \\
  \texttt{amirhosk@uci.edu} \\
   \And
 Tryphon T.\ Georgiou \\
  Department of Mechanical\\
  and Aerospace Engineering\\
  University of California, Irvine\\
  California, USA \\
  \texttt{tryphon@uci.edu} \\
}

\begin{document}
\maketitle

\begin{abstract}
This manuscript introduces a regression-type formulation for approximating the Perron-Frobenius Operator by relying on distributional snapshots of data. These snapshots may represent densities of particles. The Wasserstein metric is leveraged to define a suitable functional optimization in the space of distributions. The formulation allows seeking suitable dynamics so as to interpolate the distributional flow in function space. A first-order necessary condition for optimality is derived and utilized to construct a gradient flow approximating algorithm. The framework is exemplied with numerical simulations.
\end{abstract}


\section{Introduction}

It is often the case that dynamics are to be inferred by the collective response
of dynamical systems (particles, agents, and so on) recorded as distributional snapshots of observables~\cite{Klus_2018}. Regardless of whether the underlying dynamics is linear or not, provided there is no interaction between particles, the distributional data on observables evolve under the action of a linear operator. 
The two broadly-studied alternatives for this purpose are the Perron-Frobenius and the Koopman operators, both known as transfer operators. They are indeed linear, but defined on infinite-dimensional spaces of distributions and of observable (functions), respectively, and are adjoint to one another~\cite{klus2015numerical}.

Modeling and approximation of transfer operators often relies on samples of along collections of trajectories, e.g., see
\cite{mezic2020numerical,kutz2016dynamic,
klus2020data,berry2020bridging}.
This, in fluid mechanical systems, can be effected via recording the motion tracers seeded in the flow; such tracers provide pointwise correspondence among particles at different snapshots. However, perhaps equally often, in many real-world situations, complete trajectories may not available. Labeling and tracking particles individually is simply not feasible. In such cases, distributions of ensembles at different time instances is the only accessible data. This may also be the case
in applications, as in modeling flow/traffic, when density, average speed, and other parameters quantifying congestion are being recorded and available, and not the path of individual drivers.
Herein, we are concerned with such problems where dynamics are to be inferred from data on density flows. We advance a viewpoint that leverages the geometry of {\em Optimal Mass Transport} (OMT) and the {\em Wasserstein metric} on distributions, to identify underlying dynamics.

Besides applicactions related to flows of particles and collections of dynamical systems, the problems we consider are relevant in image registration, tumor growth monitoring, and system identification from visual data~\cite{sbalzarini2013modeling}. Another instance is domain adaptation, which aims at finding a model on a {\em target data} distribution, by training on a {\em source data}
distribution~\cite{courty2014domain,yair2019parallel}.

A popular and effective method in identifying dynamics using snapshots of data is due to Schmid and Sesterhenn (2008) and is known as DMD ({\em Dynamic Mode Decomposition}). Their algorithm aimed at modeling time-series measurements of fluid flow data~\cite{schmid2008dynamic}.  
The connection between DMD and the Koopman operator was pointed out and discussed in~\cite{rowley2009spectral};
a reformulation as a least-squares regression problem was proposed in 2014~\cite{tu2013dynamic} and
an extension, referred to as {\em extended DMD}, for approximating the eigenvalues and eigenfunctions of the Koopman operator was proposed in \cite{williams2015data,williams2014kernel}.

The Liouville operator~\cite{rosenfeld2019occupation} is another example of a linear operator associated with non-linear dynamics; this is the infinitesimal generator for the Koopman operator~\cite{rosenfeld2019dynamic}. In this context, we also mention the concept of occupation kernels which allows for the embedding of a dynamical system into a {\em Reproducing Kernel Hilbert Space} (RKHS). For further studies and taxonomy of the substantial and rapidly expanding literature we refer to~\cite{cunningham2015linear}. 

A well-known method for the approximation of \PFO$\,$ is Ulam's method, in which the evolution of a set of test points within the discretized state-space under the action of dynamics leads to a probability matrix in the discretized state-space~\cite{li1976finite,froyland2014computational}. There are other methods to approximate Perron-Frobenius operator, most of which rely on Petrov-Galerkin projections of infinite-dimensional operators onto some finite-dimensional subspace (see for example \cite{klus2015numerical,bose2006dynamical,ding1998maximum}). Also, one can utilize one of the aforementioned techniques to approximate the Koopman operator and use the duality property to find an approximate representation for the Perron-Frobenius operator~\cite{Klus_2018}.  
These approaches hypothesize the existence of pointwise correspondence among the distributions at different snapshots as the data are collected along one or several trajectories of the dynamics. 

In this paper, data are assumed to be probability distributions over a suitable state-space, and that any statistical dependence between pairs of distributions is not available. 
These observations (one-time marginal distributions) are the successive projections of the flow generated by the underlying dynamics. We seek a suitable approximation of \PFO $\,$
and, thereby, an embedding of the dynamics into a function space based on these distributional snapshots. The Wasserstein metric is employed to define an appropriate cost, by minimization of which, a desirable embedding can be achieved. 
This notion of distance, which represents cost of transport, compares two probability distributions based on the ground metric of the underlying state-space. 
The Wasserstein metric is becoming increasingly popular in recent years due to a number of natural and useful properties (e.g., being weakly continuous, allowing efficient computation via entropic regularization)~\cite{villani2008optimal,benamou2015iterative,chen2016relation}.

The paper is organized as follows. Notation and preliminaries on transfer operators are presented in Section  \ref{sec:III}, and rudiments of Wasserstein geometry needed for the development of the method are explained in Section \ref{sec:IV}. These tools are then used in Section \ref{sec:V} to derive a first-order necessary condition for two different approximations of Perron-Frobenius operator. Further, a gradient-descent approach in finding sought parameters in a system identification setting is presented. The proposed framework is highlighted via two numerical examples in Section \ref{sec:VI}.

\section{Transfer operators}\label{sec:III}

In this section we discuss the \PFO$\,$  and Koopman operators. These encode information on the underlying dynamical equations, which are nonlinear, in general. The operators are linear albeit on infinite-dimensional spaces, the space of distributions and observables, respectively.  
Although our study focuses on approximating the \PFO, we concisely summarize the duality between the two \cite{klus2015numerical}.   

\subsection{Notation}\label{sec:notation}
The three-tuple $(\X,\Sigma,\lambda)$ represents a measure space $\X\subset \Rd$ equipped with a sigma-algebra $\Sigma$ and measure $\lambda$. Typically, and unless otherwise stated, $\X=\Rd$, $\Sigma$ is the Borel algebra, and $\lambda$ the Lebesgue measure. 
The Banach space $L^p(\X)$ ($1\leq p \leq \infty$) is the space of $p$-Lebesgue integrable functions endowed with the norm $\|\cdot \|_{L^p}$. 
We denote by $(\PD,W_2)$ the Wasserstein space where $\PD$ is the set of Borel probability measures with finite second moments, and $W_2$ the Wasserstein distance.
The push-forward of a measure $\nu$ by the measurable map $S : \X \to \X$ is denoted by $\nu'=S_\#\nu \in \PD$, meaning
$\nu'(B) = \nu(S^{-1}(B))$ for every Borel set $B$.
If a measure $\mu_f\in \PD$ is absolutely continuous with respect to the Lebesgue measure, 
then we can assign to $\mu_f$, a density $f\in L^1(\X)$, that is, a positive function with unit $L^1$-norm, such that $\mu_f(B)=\int_B f d\lambda$, for every Borel set $B$.
 The Dirac measure at point $x$ is denoted by $\delta_{x}$. 

\subsection{\PFO}
A discrete-time dynamical system 
\[
x_{k+1}=S(x_k)
\]
 on $\X$ is defined by a $\lambda$-measurable state transition map
$
S: \X \to \X.
$
This map is assumed to be non-singular throughout this paper, which guarantees that the push-forward operator under $S$ preserves the absolute continuity of (probability) measures with respect to $\lambda$.  The time is assumed to be discrete. In other words, for the time lag $\tau$, the evolution of measures under $S$ can be written as $\mu_{t_k+\tau}=S_\#{\mu_{t_k}}$, ($k=1,2,\ldots$); for convenience we compress the notation by writing ${\mu_{t_k}}=:\mu_k$.

The \PFO ~(PFO), $P:L^1(\mathbb{X})\rightarrow L^1(\mathbb{X})$, is defined by
\[
\int_A Pf~ d\lambda=\int_{S^{-1}(A)} f~ d\lambda,\quad \forall A\in\Sigma
\]
for $f\in L^1(\mathbb{X})$. When $f$ is a density associated with the probability measure $\mu_f$, PFO can be thought of as a push-forward map, that is, $P\mu_f=S_\# \mu_f$.
The connection between the dynamics and PFO can be seen in that the PFO translates the center of a Dirac measure $\delta_x\in L^1(\X)$ in compliance with the underlying dynamics, that is, $S_\#\delta_x=\delta_{S(x)}$. 

It is standard that PFO is a Markov operator, namely, a linear operator which maps probability densities to probability densities. It is also a weak contraction (non-expansive map), in that, $\|Pf\|_{L^1}\leq \|f\|_{L^1}$ for any $f \in L^1(\X)$. For many dynamical systems, the PFO drives the densities into an invariant one (measure, in general) which is unique if the map $S$ is ergodic with respect to $\lambda$.

\subsection{Koopman operator}
The Koopman operator (KO) 
with respect to S, $U:L^\infty(\X)\to L^\infty(\X)$, is the infinite-dimensional linear operator
\[
Uf(x)=f(S(x)),~~ \forall x\in\X,~~\forall f\in L^\infty(\X),
\]
see e.g., \cite{budivsic2012applied}.
This is a positive operator and a weak contraction, that is, $\|Uf\|_{L^\infty}\leq \|f\|_{L^\infty}$ for any $f \in L^\infty(\X)$.

It is straightforward to see that KO is the dual of PFO, namely, \[
\langle Pf,g\rangle _\lambda=\langle f,Ug\rangle _\lambda,~~\forall f\in L^1(\X),~~g\in L^\infty(\X)
\]
where $\langle\cdot,\cdot\rangle_\lambda$ is the duality pairing between $L^1(\X)$ and $L^\infty (\X)$.
To reconstruct the underlying dynamics $(S)$ from KO, we can pick the full-state observable $g(x)=x$, where $g$ is a vector-valued observable and KO acts on it componentwise.

\subsection{Data-driven approximation of transfer operators}
As mentioned earlier, the most popular method in the literature to discretize PFO is the Ulam's method \cite{li1976finite,froyland2014computational}. In this method, the state-space ($\X$) is divided into a finite number of disjoint measurable boxes $\{B_1, . . . , B_n\}$. The PFO is approximated with a $n\times n$ matrix with elements $p_{ij}$. To do so, first we choose a large number ($k$) of test points $\{x_l^i\}_{l=1}^k$ within each Box $B_i$ randomly. Then, the elements of this matrix can be estimated by

\[
p_{ij}=\frac{1}{k}\sum_{l=1}^{k} \mathbf{1}_{B_j}(S(x_l^i))
\]
where $\mathbf{1}_{B_j}$ denotes the indicator function for the box $B_j$. 

Extended dynamic mode decomposition (EDMD)~\cite{williams2015data}, on the other hand, approximates the Koopman operator for an available time series of data, i.e., $\{x_i\}_{i=1}^m$. First, a dictionary of observables $D=\{\phi_i(\cdot)\}_{i=1}^k$ is chosen. We then consider the vector-valued function $\Phi=[\phi_1 ~\phi_2~ \ldots ~\phi_k]^T$. We stack up the values of this function at the snapshots in two matrices as \begin{align*}
    &\Phi_{[1,m-1]}=[\Phi(x_1)~\ldots~\Phi(x_{m-1})],\\ &\Phi_{[2,m]}=[\Phi(x_2)~\ldots~\Phi(x_{m})].
\end{align*} 
A finite-dimensional approximation of the restriction of the Koopman operator on the span of $D$ can be sought by considering a $k\times k$ matrix $K$ that satisfies 
\begin{equation}
    \label{EDMD}
    \Phi_{[2,m]}=K\Phi_{[1,m-1]}.
\end{equation}
Depending on the values of $m$ and $k$, the system of equations \eqref{EDMD}, may be over- or under-determined. For example, if it is over-determined, $K$ can be obtained by solving a corresponding least-squares problem.


\section{Rudiments of Wasserstein space}\label{sec:IV}
In this section, we recall the definition and some properties of the Wasserstein distance  \cite{villani2003topics,ambrosio2013user}, which are used in this paper. 

Let $\mu_0$ and $\mu_1$ be two probability measures in $\PD$. In the Monge's formulation of optimal transport, a mapping $T^*: ~\X \rightarrow \X$ is sought such that $T^*_\#\mu_0 = \mu_1$ and
\begin{equation*}
  \int_{\X} \|T^*(x)-x\|_2^2 ~d\mu_0  \leq \int_{\X} \|T(x)-x\|_2^2 ~d\mu_0
\end{equation*}
for any transport map $T$ such that $T_\#\mu_0 = \mu_1$. This is the minimization of a quadratic cost over the space of maps $T:~\X \rightarrow \X$
which ``transport'' mass $d\mu_0(x)$ at $x$ so as to match the final distribution $\mu_1$. If $\mu_0$ and $\mu_1$ are absolutely continuous, Brenier's characterization states that the optimal transport problem has a unique solution obtained as gradient of a convex function $\phi$, that is a monotone map $T^*=\nabla\phi(x)$ \cite{brenier1987decomposition}.

In case a transport map fails to exist, as is the case when $\mu_0$ is a discrete probability measure and $\mu_1$ is absolutely continuous, we consider a relaxation of Monge's problem, known as the Kantorovich's formulation, in which one seeks a joint distribution (referred to as coupling) $\pi$ on $\X\times\X$, having marginals $\mu_0$ and $\mu_1$ along the two coordinates, namely,
\begin{equation*}
W_2^2(\mu_0, \mu_1):=    \inf_{\pi \in \Pi(\mu_0, \mu_1)} \int_{\X\times\X} \|{x-y}\|^2 d\pi(x,y)
\end{equation*}
where $\Pi(\mu_0, \mu_1)$ is the space of ``couplings'' with marginals $\mu_0$ and $\mu_1$. In this, a minimizer always exists, and we use $\Pi^*(\mu_0, \mu_1)$ to denote the space of optimal couplings between the marginals $\mu_0$ and $\mu_1$.
In case the optimal transport map for the Monge problem exists, the consistency between the two problems can be realized through the relation $\pi = (x,T^*(x))_\#\mu_0$. 

The square root of the optimal cost, namely $W_2(\mu_0, \mu_1)$, defines a metric on $\PD$ referred to as the Wasserstein metric \cite{ambrosio2004gradient,villani2008optimal}. Moreover, assuming that $T^*$ exists, the constant-speed geodesic between $\mu_0$ and $\mu_1$ is given by
\begin{equation*}
    \mu_t=\{(1-t)x+tT^*(x)\}_\#\mu_0,~~0\leq t \leq 1,
\end{equation*}
and known as {\em McCann's  displacement interpolation} \cite{mccann1997convexity}.

In the following, we state an important lemma from measure theory which will be used in the proof of main theorem in this paper.
\begin{lemma}[Gluing lemma~\cite{ambrosio2004gradient,villani2008optimal}]
\label{gluing}
Let $\X_1$, $\X_2$, and $\X_3$ be three copies of $\X$. 
Given three probability measures $\mu_i(x_i)\in \mathcal{P}_2(\mathbb{X}_i),~i=1,2,3$ and the couplings $\pi_{12}\in \Pi(\mu_1,\mu_2)$, and $\pi_{13}\in\Pi(\mu_1,\mu_3)$, there exists a probability measure $\pi(x_1,x_2,x_3)\in \mathcal{P}_2(\X_1\times \X_2 \times \X_3)$ such that $(x_1,x_2)_\#\pi=\pi_{12}$ and $(x_1,x_3)_\#\pi=\pi_{13}$. Furthermore, the measure $\pi$ is unique if either $\pi_{12}$ or $\pi_{13}$ are induced by a transport map.
\end{lemma}
That is, the gluing lemma states that for any two given couplings, which are consistent along one coordinate, we can find a measure on the product space $(\X_1\times \X_2 \times \X_3)$ whose projections onto each pair of coordinates match the given couplings, respectively. 
With this, we are ready to present the main results in the next section.

\section{Main results}\label{sec:V}
In this section,  we formally define the problem of PFO approximation in the presence of distributional snapshots for a dynamical system. 
As already noted, it is assumed that there is no information on the correlation between each pair of data points (distributions). 
We seek system dynamics, $S:\X\to \X$, as a $\lambda$-measurable map 
such that it can serve as a model for the flow encoded in the sequence of data points $\mu_1$, $\mu_2$, \ldots, $\mu_m$. This is in the sense that, either $S_\# \mu_k=\mu_{k+1}$ over the data set for $k\in\{1,\ldots,m-1\}$ (exact matching), or that the discrepancy between $S_\# \mu_k$ and $\mu_{k+1}$, for the successive data points, is small in the average over the available record of distributions.
Below, in Section \ref{sec:one}, we first develop the case where $S$ is a linear map 
\[
S\;:\;x\mapsto Ax,
\]
 with $A\in \mathbb R^{d\times d}$. Then, in Section \ref{sec:two}, we detail the approach for the case where $S(\cdot) = \sum_{j=1}^{n}\theta_j y_j(\cdot)$ is nonlinear (in general) expressed in terms of a linear combination of specified basis functions $y_j$, $j\in\{1,\ldots,n\}$.

\subsection{First-order approximation}\label{sec:one}

We first draw an analogy with the EDMD problem by stating the problem to find a matrix that satisfies the condition in Eq. \eqref{EDMD}.
Thus, given a sequence of probability measures $\{\mu_i\}_{i=1}^m$ in $\PD$, we seek to find a matrix $A\in M(d)$ (the space of real $d\times d$ matrices) such that 
\begin{equation}
    \label{WDMD}
    [\mu_2~\mu_3~\ldots~\mu_m]=(Ax)_\#  [\mu_1~\mu_2~\ldots~\mu_{m-1}].
\end{equation}
In \eqref{WDMD}, similar to EDMD, the probability distributions ($\mu_1$, $\mu_2$,\ldots) are stacked in arrays, where one is the shifted version of the other. The push-forward operator acts on ``stacked up'' measures separately. 

Typically, the problem is over-determined, in which case there might not exist a matrix $A$ that satisfies \eqref{WDMD}, we consider the following regression-type formulation.
\begin{problem}\label{prob2}
Determine a matrix $A\in M(d)$ that minimizes
\begin{equation} 
\label{W_reg}
F(A)= \sum_{i=1}^{m-1} W^2_2(Ax_{\#}\mu_{i},\mu_{i+1}).
\end{equation}
\end{problem}

If \eqref{WDMD} has a solution, it trivially coincides with the minimizer of Problem \ref{prob2} and $F(A)=0$. 

If, on the other hand, all the measures are Dirac, that is, $\mu_i=\delta_{x_i},~i=1,\ldots,m$,
the problem to satisfy \eqref{WDMD} reduces to an ordinary DMD problem. This shows the consistency of DMD with our formulation on measures.


Next, we provide a stationarity condition that can be used to obtain the solution to Problem \ref{prob2}.

\begin{theorem}
\label{theorem1}
Consider a sequence of absolutely continuous probability measures $\{\mu_i\}_{i=1}^m$ in $\PD$. If a minimizer $A\in M(d)$ for \eqref{W_reg} exists and is nonsingular, then there exist unique
$\eta_i(x_i,x_{i+1})\in \Pi(\mu_i,\mu_{i+1})$ for each $i\in\{1,\ldots,m\}$ such that 
\[(Ax_i,x_{i+1})_{\#}\eta_i \in \Pi^*({Ax_i}_{\#}\mu_i,\mu_{i+1}),\]
and moreover, $A$ satisfies
\begin{equation}
    \sum_{i=1}^{m-1}\int_{\X\times \X} (Ax_{i}-x_{i+1})x_i^T d\eta_i(x_i,x_{i+1})=0.
\end{equation}
\end{theorem}


In the theorem, each probability measure $\eta_i$ is a coupling between two distributional snapshots $\mu_i$ and $\mu_{i+1}$ such that the push-forward measure $(Ax_i,x_{i+1})_{\#}\eta_i$ is an optimal coupling between its marginals. In turn, since these marginals are absolutely continuous by virtue of the fact that $A$ is nonsingular, the latter coupling (i.e., $(Ax_i,x_{i+1})_{\#}\eta_i$) is singular and ``sits'' on the graph of a ``Monge map.'' As explained in the proof of the theorem, application of the Gluing lemma 
shows that each $\eta_i$ exists and is unique. At this point, the absolute continuity of the marginals is essential; later on, we will discuss how to relax this assumption so as to include a class of discrete measures as well.

\noindent
{\em Proof of Theorem \ref{theorem1}}:
According to the assumption that $A$ is a minimizer of \eqref{W_reg}, the Fermat's condition 
\begin{equation}
    \label{Fermat's}
    \frac{d}{d\epsilon} F(A+\epsilon\delta A)|_{\epsilon=0}=0
\end{equation}
holds for any tangent direction $\delta A$, that is, any matrix in $M(d)$. Without loss of generality, we consider only one of the terms in \eqref{W_reg} and define 
\[
G(A)=W^2_2(Ax_{\#}\mu_1,\mu_2).
\]
To calculate the directional derivative (Gateaux derivative) of $G(A)$, first we show that for any real $\epsilon$ and $\delta A\in M(d)$
\begin{align}
\label{limsup}
&{G(A+\epsilon\delta A)-G(A)}\leq \\ \nonumber
&\langle \int_{\X\times \X} 2(Ax_{1}-x_{2})x_1^T d\eta_1(x_1,x_{2}),\epsilon \delta A\rangle _F+O(\epsilon^2)
\end{align}
where $\langle \cdot,\cdot \rangle _F$ is the Frobenius inner product and $\eta_1$ is as stated in the theorem. To do so, let the measure $\gamma_1(x_1,x_1',x_2)\in \mathcal{P}_2(\X^3)$ be such that $(x_1,x_1')_\#\gamma_1=(x_1,Ax_1)_\#\mu_1$ and $(x_1',x_2)_\#\gamma_1\in \Pi^*({Ax_1}_\#\mu_1,\mu_2)$. Since these two constraints coincide along $x_1'$, by application of the Gluing lemma, we conclude that $\gamma_1$ exists. Moreover, as the projection of $\gamma_1$ onto $(x_1',x_2)$ is the optimal coupling between two absolutely continuous measures, it is induced by a transport map (Monge map), and thus the choice of $\gamma_1$ is unique by once again invoking the Gluing lemma. Then, $\eta_1:=(x_1,x_2)_\#\gamma_1$ where its uniqueness immediately results from that of $\gamma_1$. Hence,
\begin{align*}
    &G(A+\epsilon\delta A)-G(A)\leq \\
    &\int_{\X_1\times\X_{2}}  (\|(A+\epsilon \delta A)x_1-x_2\|_2^2-\|Ax_1-x_2\|_2^2)d\eta_1(x_1,x_{2}).
\end{align*}
This follows from the fact that $G(A+\epsilon \delta A)$ is the Wasserstein distance (i.e., the minimum among all the couplings between ${(A+\epsilon \delta A)x_1}_\#\mu_1$ and $\mu_2$). Finally, by expanding the integrand  above with respect to $\epsilon$, \eqref{limsup} is derived. 

Without loss of generality we take $\epsilon>0$. According to \eqref{limsup}, we can readily conclude that
\begin{align*}
&\limsup\limits_{\epsilon\to 0} \frac{G(A+\epsilon\delta A)-G(A)}{\epsilon}\leq \\ \nonumber
&\langle \int_{\X_1\times \X_{2}} 2(Ax_{1}-x_{2})x_1^T d\eta_1(x_1,x_{2}),\delta A\rangle _F.
\end{align*}

The next step of proof is to show that 
\begin{align*}
&\liminf\limits_{\epsilon\to 0} \frac{G(A+\epsilon\delta A)-G(A)}{\epsilon}\geq \\ \nonumber
&\langle \int_{\X_1\times \X_{2}} 2(Ax_{1}-x_{2})x_1^T d\eta_1(x_1,x_{2}),\delta A\rangle _F.
\end{align*}
This last inequality follows from the semi-concavity of the squared Wasserstein distance   \cite[Proposition 7.3.6]{ambrosio2008gradient}. 

By combining the ``$\liminf$" and ``$\limsup$" results, it readily follows that
\begin{align}
\label{grad_type}
&\frac{d}{d\epsilon} G(A+\epsilon\delta A)|_{\epsilon=0}=\\ \nonumber \nonumber
&\langle \int_{\X\times \X} 2(Ax_{1}-x_{2})x_1^T d\eta_1(x_1,x_{2}),\delta A\rangle _F.
\end{align}
Finally, writing the directional derivative for all the terms in \eqref{W_reg} and using Fermat's condition the proof is complete. 
\hfill $\Box$

\begin{remark}{\rm
In the statement of Theorem \ref{theorem1} we assume the existence of a minimizer $A$ to Problem \ref{prob2}. We now explain that this assumption holds in many reasonable settings, as for instance, in the case where the probability measures have compact support. To see this, note that $F(A)$ is coercive, i.e., $F(A)\to +\infty$ as $\|A\|_F\to +\infty$ for absolutely continuous $\mu_i$'s with compact support. Further, using the lower semi-continuity of $W_2$ (see Proposition 7.1.3 and Lemma 5.2.1 in \cite{ambrosio2008gradient}), we conclude the lower semi-continuity of $F(A)$ with respect to the Frobenius norm. These two observations guarantee the existence of a solution to Problem \ref{prob2}. }\hfill $\Box$
\end{remark}

\begin{remark}{\rm 
Equation \eqref{grad_type} shows how to generate a gradient flow, and thereby a steepest descent direction for minimizing $F(A)$. Specifically,
\begin{equation}
\label{gradeint_descent}
    \nabla_AF(A)=2\sum_{i=1}^{m-1}\int_{\X_i\times \X_{i+1}} (Ax_{i}-x_{i+1})x_i^T d\eta_i(x_i,x_{i+1}),
\end{equation}
allows us to construct a gradient-type numerical optimization to find the minimizer of~\eqref{W_reg}.}\hfill $\Box$
\end{remark} 

\begin{remark}{\rm
We note in passing that the setting of our approximation Problem \ref{prob2}, can be used to construct pseudo-metrics for various applications. Specifically,
an admissible set of transformations $\mathcal F$ may be available (e.g., rotations, translations, scalings of images and so on), and that these are natural for the problem at hand, and thought to ``incur no cost.'' Thence,
a distance can be defined between distributions as follows
\[
W^2_{\mathcal{F}}(\mu_0,\mu_1)=\inf_{S\in\mathcal{F}} W^2_2(S_\#{\mu_0},\mu_1).
\]
Such a construction is relevant in image registration where alignment/scaling may be desired.}\hfill $\Box$
\end{remark}

\subsection{Higher-order approximations}\label{sec:two}
In this subsection, we extend the previous result to non-linear models for the underlying dynamics.

We consider system dynamics, $S:\X\to \X$, a $\lambda$-measurable map, to be expressed as a linear combination of basis functions $y_j : \X \to \X$, with $j\in\{1,\ldots,n\}$, i.e.,
\[
S(x;\Theta)=\sum_{j=1}^n\theta_j y_j(x).
\]
where $\Theta=[\theta_1~ \ldots~\theta_n]^T\in \mathbb{R}^n$.

The set of basis functions may be chosen to include polynomials. In such a case,  the corresponding-order moments of the distributional snapshots need to exist, so that integrals remain finite. 

Extending \eqref{W_reg} to this new setting, we now consider the problem to minimize 
\begin{equation} 
\label{W_reg_F}
F(\Theta)= \sum_{i=1}^{m-1} W^2_2(S(x;\Theta)_{\#}\mu_{i},\mu_{i+1}),
\end{equation}
over $\Theta\in \mathbb{R}^n$.
We follow a strategy that is similar to that in the proof of Theorem \ref{theorem1}, to derive a first-order optimality condition for $\Theta$ in the form
\begin{equation}
    \label{higher-approx}
     \sum_{i=1}^{m-1}\int_{\X\times \X} (Y(x_i))^T(S(x_i;\Theta)-x_{i+1}) d\eta_i(x_ix_{i+1})=0.
\end{equation}
Here, $Y(x_i)=[y_1(x_i)~\ldots~y_n(x_i)]\in \mathbb{R}^{d\times n}$ and, as before, $\eta_i(x_i,x_{i+1})\in \Pi(\mu_i,\mu_{i+1})$ is such that 
\[(S(x_i;\Theta),x_{i+1})_{\#}\eta_i \in \Pi^*({S(x_i;\Theta)}_{\#}\mu_i,\mu_{i+1}).\] 
In a similar manner, the absolute continuity of $\mu_i$'s guarantees the existence and uniqueness of all the $\eta_i$'s.

Equation \eqref{higher-approx} extends our formalism to nonlinear dynamics, parametrized by the span of $Y$, for approximating the PFO. In a way similar to \eqref{gradeint_descent}, we consider the gradient of $F(\Theta)$ in \eqref{W_reg_F} with respect to $\Theta$,
\begin{equation}
\label{gradeint_descent_nonlinear}
    \nabla_{\Theta}F=2\sum_{i=1}^{m-1}\int_{\X\times \X}  (Y(x_i))^T(S(x_i;\Theta)-x_{i+1}) d\eta_i(x_i,x_{i+1}),
\end{equation}
and employ a gradient-type descent to find the minimizing value for $\Theta$.

\section{Simulation results}\label{sec:VI}
\subsection{Gaussian distributions}

\begin{figure}[htb]
	\centering
	\includegraphics[width=7cm]{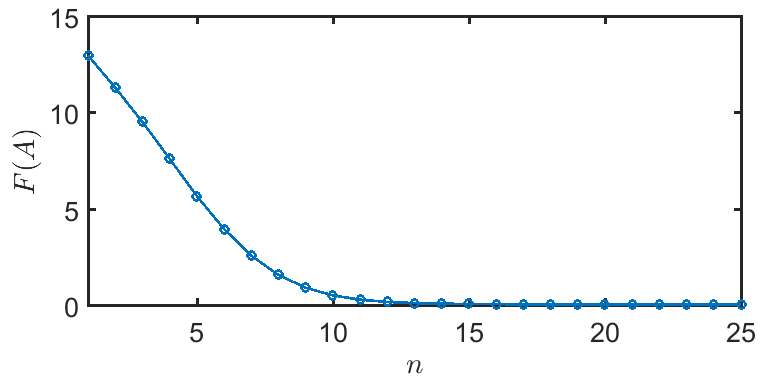}
	\caption{Value $F(A_n)$ as a function of iterated steps in \eqref{gradient_algo}.}
	\label{cost_iter}
\end{figure}

\begin{figure}[htb]
	\centering
	\includegraphics[width=\textwidth,height=4cm]{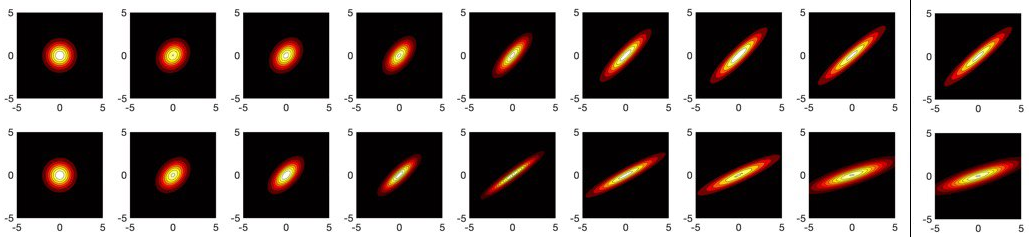}
	\caption{The rows exemplify the convergence of 
$(A_nx)_{\#}\mu_1\to \mu_2$ and $(A_nx)_{\#}\left((A_nx)_{\#}\mu_1\right)\to \mu_3$, respectively, as $n=1,\ldots,8$, towards $\mu_2$ and $\mu_3$, which are displayed on the right and separated by a vertical line (with $\mu_2$ on top of $\mu_3$).}
	\label{convg}
\end{figure}

We exemplify our framework with numerical results for the case where the distributional snapshots are Gaussian. In this case, the Wasserstein distance between distributions can be written in closed-form. 

Consider\footnote{$\mathcal N(m,C)$ denotes a Gaussian distribution with mean $m$ and covariance $C$} $\mu_0 = \mathcal N({m}_0,{C}_0)$ and $\mu_1 = \mathcal N({m}_1,{C}_1)$.
The transportation problem admits a solution in closed-form~\cite{malago2018wasserstein,masarotto2019procrustes}, with transportation (Monge) map
\begin{equation*}T^*:x\to  C_0^{-1}(C_0C_1)^{1/2}={{C}_0^{-1/2}({C}_0^{1/2}{C}_1{C}_0^{1/2})^{1/2}{C}_0^{-1/2}x},
\end{equation*}
and transportation cost $W_2(\mu_0,\mu_1)$ given by
\begin{equation}
\label{Gaussian_W2}
    { \sqrt{\|{m}_0-{m}_1\|^2+{\rm tr}({C}_0+{C}_1-2 (C_1^{1/2}C_0C_1^{1/2})^{1/2})} }
\end{equation}
where ${\rm tr}(.)$ stands for trace.

We begin with a collection $\mu_i = \mathcal N(0,C_i),~i=1,\ldots,m$ as our distributional snapshots; for simplicity we have assumed zero-means.
The cost \eqref{W_reg} reads
\begin{equation}
    \label{Gaussian_cost}
    F(A)=\sum_{i=1}^{m-1} tr(AC_iA^T+C_{i+1}-2 (C_{i+1}^{1/2}AC_iA^T C_{i+1}^{1/2})^{1/2}).
\end{equation}
The gradient $\nabla_A F(A)$, for the case of Gaussian snapshots, is expressed below directly in terms of the data $C_i$, $i\in\{1,\ldots,m\}$.
\begin{proposition}
Given Gaussian distributions $\mu_i = \mathcal N(0,C_i),~i=1,\ldots,m$, and a non-singular $A\in M(d)$,
\begin{equation}
    \label{Gaussian_gradient}
   \nabla_A F=2 \left\{A\sum_{i=1}^{m-1} C_i-(\sum_{i=1}^{m-1}(C_{i+1}AC_iA^T)^{1/2})A^{-T}\right\}.
\end{equation}
\end{proposition}

To determine a minimizer for \eqref{Gaussian_cost}, we utilize a first-order iterative algorithm, taking steps proportional to the negative of the gradient in \eqref{Gaussian_gradient}, namely,
\begin{equation}
\label{gradient_algo}
A_{n+1}=A_n-\alpha \nabla_{A} F(A_n),~~n=1,2,\ldots    
\end{equation}
for a small learning rate $\alpha>0$.

As a guiding example, and for the sake of visualization, we consider the two-dimensional state-space $\X=\R^2$, in which probability measures are evolving according to linear non-deterministic dynamics, 
\[
x_{k+1}=\left[\begin{array}{cc}
     -\frac12 & 2\\
     -1 & \frac{3}{2}
\end{array} \right]x_k+\frac{2}{5}
\left[\begin{array}{c}
     \Delta \omega_k^1  \\
     \Delta \omega_k^2 
\end{array} \right],~~k=1,2,\ldots
\]
starting from $\mu_1= \mathcal N(0,I_2)$, with $I_2$ a $2\times 2$ identity matrix. We take $\Delta \omega_k^1,~\Delta \omega_k^2 = \mathcal N(0,1)$ to be independent white noise processes. 

This dynamical system is an example of a first-order autoregressive process (AR(1)) which can also be thought of as an Euler-Maruyama approximation of a two-dimensional Ornstein-Uhlenbeck stochastic differential equation where $\Delta \omega_k^1$ and $\Delta \omega_k^2$ are the increments of two independent Wiener processes with unit step size.

We note that $A$ is neither symmetric nor positive definite, which implies that it is not a ``Monge map'' and, thus, the flow of distributions is not a geodesic path in the Wasserstein metric. 

Using the first five iterates $(m=6)$, we employ \eqref{gradient_algo} to obtain dynamics
solely on the basis of these 5 distributional snapshots. We initialize \eqref{gradient_algo} with $\alpha=0.1$ and experimented with different starting choices for $A_1$. Specifically,
we took $A_1$ to be the identity matrix $I_2$, and also, 
 the average $ A_1=\frac{1}{m-1}\sum_{i=1}^{m-1} C_i^{-1}(C_iC_{i+1})^{1/2}$, without any perceptible difference in the convergence to a minimizer. For the first choice, $A_1=I_2$, the values of $F(A_n)$ in successive iterations is shown in Fig.~\ref{cost_iter}.
%

Our data $C_i$ ($i\in\{1,\ldots,6\}$) is generated starting from $\mu_1 = \mathcal N(0,C_1)$ with $C_1=I_2$, i.e., the $2\times 2$ identity matrix, and the gradient search for the minimizer is initialized using $A_1=I_2$ as well. 
In Fig.\ \ref{convg} we display contours of probability distributions.
Specifically, on the right hand side, separated by a vertical line, we display the contours for $\mu_2 = \mathcal N(0,C_2)$
and $\mu_3 = \mathcal N(0,C_3)$, with $\mu_2$ on top of $\mu_3$.
Then, horizontally, from left to right, we display contours corresponding to the
approximating sequence of distributions. The first row exemplifies the convergence
\[
(A_nx)_{\#}\mu_1\to \mu_2,
\]
whereas the second row, exemplifies the convergence
\[
(A_nx)_{\#}\left((A_nx)_{\#}\mu_1\right)\to \mu_3,
\]
as $n=1,\ldots,8$.



\begin{figure}[htb]  
\centering
\subfigure[]{\resizebox{!}{3.2cm}{\includegraphics{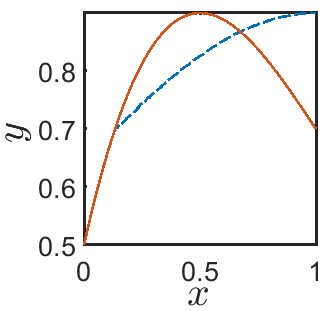}}}
\subfigure[]{\resizebox{!}{3.2cm}{\includegraphics{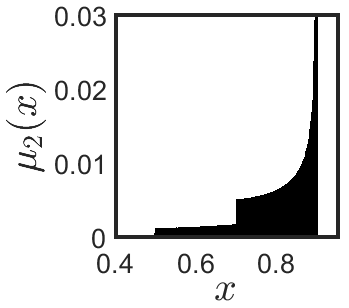}}}
\caption{The two maps in (a) transport a uniform distribution on $[0,1]$ to the same discontinuous density in (b). Monge map (blue) is injective but not in $C^1$ everywhere. The non-injective map (red) is in $C^1$.}
\label{Monge_map}
\end{figure}

\begin{figure}[htb]
	\centering
	\includegraphics[width=\textwidth,height=3cm]{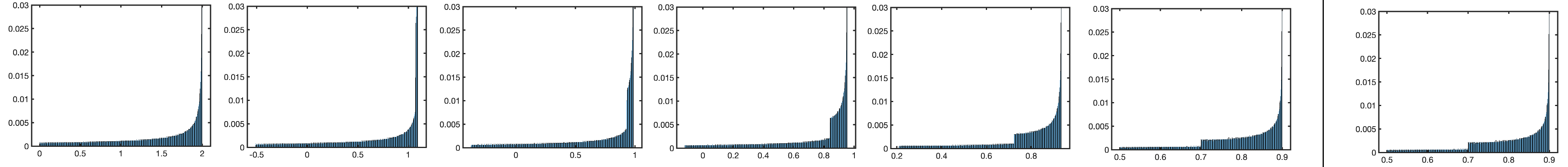}
	\caption{ 
	The evolution of uniform distribution under $S(x;\Theta)$ at different iterations of the algorithm. On the right-hand side the target density is depicted. In the beginning (left) no jump discontinuity is observed.}
	\label{sing_1}
\end{figure}

\begin{figure*}[htb]
	\centering
	\includegraphics[width=\textwidth,height=3cm]{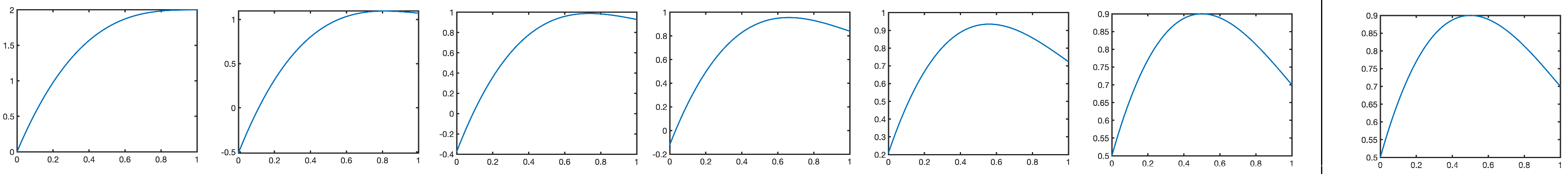}
	\caption{ 
	The transport map $S(x;\Theta)$ at different iterations of the algorithm. This shows the convergence to the non-injective map.}
	\label{sing_2}
\end{figure*}

\subsection{Non-linear dynamics}
For our second example, to highlight the use of the approach, we consider the $C^1$ (continuously differentiable) map $S\,:\,x\mapsto S(x)$, on $\X=\mathbb{R}$ with
\begin{equation}\label{eq:Sx}
S(x)=0.7+0.6 (1-x)-0.8(1-x)^3.
\end{equation}
The idea for this example has been borrowed from \cite{moosmuller2020geometric}.
The map $S$ is depicted in Fig.\ \ref{Monge_map}(a), in red solid curve, 
and pushes forward a uniform distribution on $[0,1]$ to distribution with discontinuous density. This density is shown in Fig.\ \ref{Monge_map}(b). Due to the fact that the density is discontinuous, the optimal transport (Monge) map has a ``corner'' (not smooth) and is displayed in Fig.\ \ref{Monge_map}(a), with a dashed blue curve. 

The method outlined in this paper allows us to seek a transportation map, within a suitably parametrized class of functions, that pushes forward $\mu_1$  (here, this is the uniform distribution on $[0,1]$) to $\mu_2$ displayed in Fig.\ \ref{Monge_map}(b).
To this end, we select the representation
\[
S(x;\Theta)=\theta_3+\theta_2 (1-x)+\theta_1(1-x)^3,
\]
in the basis $Y=\{1, (1-x), (1-x)^3 \}$, and seek to determine the parameters $\theta_k$ ($k\in\{1,2,3\}$) via a
gradient-descent as in \eqref{gradeint_descent_nonlinear}.

The two probability distributions are approximated using 100 sample points (drawn independently).
We initialize with $\theta_1=-2$, $\theta_2=0$, and $\theta_3=2$. A discrete optimal transport problem is solved to find the joint distributions $\eta_i$ in \eqref{gradeint_descent_nonlinear} at each time step.
The convergence is depicted in Fig.\  \ref{sing_2}, where successive iterants are displayed from left to right below the resulting pushforward distribution. On the right hand side, separated by vertical lines, the 
target $\mu_1$ is displayed above the cubic map in \eqref{eq:Sx}.

It is worth observing that,
as illustrated in Fig.\ \ref{sing_2}, our initialization corresponds to an injective map resulting in no discontinuity in the first pushforward distribution. In successive steps however, as the distributions converge to $\mu_1$ and the maps to $S(x)$ in \eqref{eq:Sx}, a discontinuity appears tied to the non-injectivity of the maps with updated parameters.

\section{Concluding remarks}\label{sec:conclusions}
We presented an approach to interpolate distributional snapshots by identifying suitable underlying dynamics. It is assumed that no information on statistical dependence between successive pairs of distributions is available. The scheme we propose aims at modeling a \PFO $\;$associated with underlying unknown dynamics. It is based on formulating a regression-type optimization problem in the Wasserstein metric, weighing in distances between successive distributional snapshots. A first-order necessary condition is derived that leads to a gradient-descent algorithm. The method extends to search for nonlinear dynamics assuming a suitable parametrization of the nonlinear state transition map in terms of selected basis functions.  Two academic examples are presented to highlight the approach as applied in two cases, the first specializing to Gaussian distributions and the second dealing with more general distributions (albeit with one-dimensional support for simplicity).

\bibliographystyle{plain}        
\bibliography{Reference}

\end{document}